\documentclass[12pt,reqno]{amsart}

\usepackage{latexsym,mathrsfs,color}

\usepackage[colorlinks=true,
linkcolor=webblue, filecolor=webbrown, citecolor=webred ]{hyperref}
\definecolor{webblue}{rgb}{0,.5,0}
\definecolor{webred}{rgb}{0,.5,0}
\definecolor{webbrown}{rgb}{.6,0,0}

\newtheorem{theorem}{Theorem}

\newtheorem{corollary}[theorem]{Corollary}
\newtheorem{proposition}[theorem]{Proposition}
\newtheorem{conjecture}[theorem]{Conjecture}

\newcommand{\mo}{\mathcal{O}}
\newcommand{\msn}{{\mathcal S}_n}
\newcommand{\exc}{{\rm exc\,}}
\newcommand{\bn}{{\rm\bf n}}
\newcommand{\rz}{{\rm RZ}}
\newcommand{\sep}{\preceq}
\newcommand{\sn}{S_n}
\newcommand{\lrf}[1]{\lfloor #1\rfloor}
\newcommand{\lrc}[1]{\lceil #1\rceil}

\begin{document}

\title{$q$-Eulerian polynomials and polynomials with only real zeros}
\author{Shi-Mei Ma}
\address{Department of Applied Mathematics,
         Dalian University of Technology,
         Dalian 116024, P. R. China}
\email{simons\_ma@yahoo.com.cn}
\author{Yi Wang}
\address{Department of Applied Mathematics,
         Dalian University of Technology,
         Dalian 116024, P. R. China}
\email{wangyi@dlut.edu.cn}
\thanks{The second author is responsible for all the communications} 
\date{\today}

\begin{abstract}
Let $f$ and $F$ be two polynomials satisfying
$F(x)=u(x)f(x)+v(x)f'(x)$. We characterize the relation between the
location and multiplicity of the real zeros of $f$ and $F$, which
generalizes and unifies many known results, including the results of
Brenti and Br\"and\'en about the $q$-Eulerian polynomials.
\end{abstract}

\maketitle

\section{Introduction}
Let $\sn$ denote the permutation group on the set $\{1,2,\ldots,n\}$
and $\pi=a_1a_2\cdots a_n\in\sn$. An {\it excedance} in $\pi$ is an
index $i$ such that $a_i>i$. Let $\exc(\pi)$ denote the number of
excedances in $\pi$. The classical Eulerian polynomials $A_n(x)$ are
defined by
$$A_0(x)=1,\quad A_n(x)=\sum_{\pi\in\sn}x^{\exc(\pi)+1}\quad\textrm{for $n\ge 1$},$$
and have been extensively investigated. It is well known that the
classical Eulerian polynomials satisfy the recurrence relation
$$A_{n+1}(x)=(n+1)xA_{n}(x)+x(1-x)A'_{n}(x)$$
(see B\'ona~\cite[p. 23]{Bona04} for instance). In~\cite{FS70},
Foata and Sch\"utzenberger introduced a $q$-analog of the classical
Eulerian polynomials defined by
$$A_0(x;q)=1,\quad A_n(x;q)=\sum_{\pi\in\msn}x^{\exc(\pi)}q^{c(\pi)}\quad\textrm{for $n\ge 1$},$$
where $c(\pi)$ is the number of cycles in $\pi$. The first few of
the $q$-Eulerian polynomials are
$$A_0(x;q)=1, A_1(x;q)=q, A_2(x;q)=q(x+q),
A_3(x;q)=q[x^2+(3q+1)x+q^2].$$ Clearly,
$A_n(x)=xA_n(x;1)$ for $n\ge 1$. 
Brenti~\cite{Bre94,Bre00} further studied $q$-Eulerian polynomials
and established the link with $q$-symmetric functions arising from
plethysm. He obtained the recurrence relation
\begin{equation}\label{anxq-rr}
A_{n+1}(x;q)=(nx+q)A_{n}(x;q)+x(1-x)\frac{d}{d x}A_{n}(x;q)
\end{equation}
(\cite[Proposition 7.2]{Bre00}) and showed that $A_n(x;q)$ has only
real nonnegative simple zeros when $q$ is a positive rational number
(\cite[Theorem 7.5]{Bre00}). He also proposed the following.
\begin{conjecture}[{\cite[Conjecture 8.8]{Bre00}}]
\label{brenti-conj} Let $n,t\in\mathbb{N}$. Then $A_n(x;-t)$ has
only real zeros.
\end{conjecture}

The conjecture has been settled recently by
Br\"and\'en~\cite{Bra06}. Let
$$E_n(x;q)=(1+x)^{n}A_n\left(\frac{x}{1+x};q\right).$$ Then it is clear that
$A_n(x;q)$ has only real zeros if and only if $E_n(x;q)$ does. The
recurrence (\ref{anxq-rr}) induces
$$E_{n+1}(x;q)=q(1+x)E_n(x;q)+x(1+x)\frac{d}{d x}E_n(x;q),$$
with $E_0(x;q)=1$.  Using multiplier $n$-sequences, Br\"and\'en can
manage to prove that if $q>0, n+q\le 0$ or $q\in\mathbb{Z}$, then
$E_n(x;q)$ has only real zeros, and so does $A_n(x;q)$ (see
\cite[Theorem 6.3]{Bra06} for details). In the next section, we will
obtain a more precise result directly by the recurrence
(\ref{anxq-rr}) as an application of our main results in this paper.

Polynomials with only real zeros arise often in combinatorics,
algebra, analysis, geometry, probability and statistics. For
example, let $S(n,k)$ be the Stirling numbers of the second kind and
$B_n(x)=\sum_{k=0}^{n}S(n,k)x^k$ the Bell polynomials. Then
\begin{equation}\label{r-bnx}
    B_n(x)=xB_{n-1}(x)+xB'_{n-1}(x),\qquad B_0(x)=1.
\end{equation}
For showing that the Stirling behavior is asymptotically normal,
Harper~\cite{Har67} showed that the Bell polynomials have only real
simple zeros by means of the recurrence (\ref{r-bnx}).

Let $\rz$ denote the set of real polynomials with only real zeros.
Furthermore, denote by $\rz(I)$ the set of such polynomials all
whose zeros are in the interval $I$. Suppose that $f,F\in\rz$. Let
$\{r_i\}$ and $\{s_j\}$ be all zeros of $f$ and $F$ in nonincreasing
order respectively. We say that $f$ {\it separates} $F$, denoted by
$f\sep F$, if $\deg f\le\deg F\le\deg f+1$ and
\begin{equation*}
s_1\ge r_1\ge s_2\ge r_2\ge s_3\ge r_3\ge\cdots.
\end{equation*}
It is well known that if $f\in\rz$, then $f'\in\rz$ and $f'\sep f$.
Following Wagner~\cite{Wag91}, a real polynomial is called {\it
standard} if it has positive leading coefficient.

Let $f$ and $F$ be two polynomials satisfying the relation
\begin{equation}\label{Fff'}
F(x)=u(x)f(x)+v(x)f'(x).
\end{equation}
A natural question is in which cases $f$ has only real zeros implies
that $F$ does. There have been some partial results
\cite{LWaam,WYjcta05}. However, these results can not tell us the
relation of the multiplicity and location of zeros of $f$ and $F$.
The main object of this paper is to provide characterizes for such a
problem, which can give a unified explanation of many known results.
\section{Main results}
In this section we present the main results of this paper.
\begin{theorem}\label{m-thm}
Let $f,F$ be two standard polynomials and satisfy
\begin{equation}\label{Fff'}
F(x)=u(x)f(x)+v(x)f'(x),
\end{equation}
where $u(x),v(x)$ are real polynomials and $\deg F=\deg f$ or $\deg
f+1$. Assume that $f\in\rz$ and $v(r)\le 0$ whenever $f(r)=0$. Then
$F\in\rz$ and $f\sep F$. Moreover, if $r$ is a zero of $f$ with the
multiplicity $m$, then the multiplicity of $r$ as a zero of $F$ is
\begin{enumerate}
  \item [(a)] $m-1$ if $v(r)\neq 0$; or
  \item [(b)] $m$ if $v(r)=0$ but $u(r)+mv'(r)\neq 0$; or
  \item [(c)] $m+1$ if $v(r)=0$ and $u(r)+mv'(r)=0$.
\end{enumerate}
Furthermore, we have the following result.
\begin{enumerate}
\item [(A)]
  Suppose that $f\in\rz(-\infty,r]$, where $r$ is the largest
  zero of $f$, with the multiplicity $m$.
  Then $F\in\rz(-\infty,r]$ if and only if $v(r)=0$ and $u(r)+mv'(r)\ge 0$.
\item [(B)]
  Suppose that $f\in\rz[r,+\infty)$, where $r$ is the
  smallest zero of $f$, with the multiplicity $m$.
  Then $F\in\rz[r,+\infty)$ if and only if $\deg F=\deg f$, or $v(r)=0$ and $u(r)+mv'(r)\le
  0$.
\end{enumerate}
\end{theorem}
\begin{proof}
The first part of the statement about $F\in\rz$ and $f\sep F$ can be
followed from \cite[Theorem 2.1]{LWaam}. However, we give a direct
proof of it for our purpose. Without loss of generality, we may
assume that $f$ and $F$ are monic. Let
$f(x)=\prod_{i=1}^{k}(x-r_i)^{m_i}$ where $r_1,\ldots,r_k$ are
distinct zeros of $f(x)$ with the multiplicities $m_1,\ldots,m_k$
respectively. Then $\prod_{i=1}^{k}(x-r_i)^{m_i-1}|F(x)$. Denote
$g(x)=\prod_{i=1}^{k}(x-r_i)$ and
$G(x)=F(x)/\prod_{i=1}^{k}(x-r_i)^{m_i-1}$. Then $\deg G-\deg g=\deg
F-\deg f=0$ or $1$, and by (\ref{Fff'}),
\begin{equation}\label{Gg}
G(x)=u(x)g(x)+v(x)\sum_{i=1}^{k}\frac{m_ig(x)}{x-r_i}.
\end{equation}

Consider first the case $v(r_i)<0$ for all $i$. Let
$r_k<\cdots<r_1$. Then by (\ref{Gg}), the sign of $G(r_i)$ is
$(-1)^i$ for $i=1,\ldots,k$. Note that $G(x)$ is monic and $\deg
G-\deg g=\deg F-\deg f$. Hence $G(x)$ has precisely one zero in each
of $k$ intervals $(r_k,r_{k-1}),\ldots,(r_2,r_1),(r_1,+\infty)$ and
has an additional zero in the interval $(-\infty,r_k)$ if $\deg
G-\deg g=1$. Thus $G\in\rz$ and $g\sep G$. It implies that $F\in\rz$
and $f\sep F$. Clearly, $r_i$ is not a zero of $G$. So $r_i$ is a
zero of $F$ with the multiplicity $m_i-1$. This proves (a).

Next consider the general case. Let $v_j(x)=v(x)-1/j$ and
$F_j(x)=u(x)f(x)+v_j(x)f'(x)$. Then $v_j(r_i)<0$ for all $i$ when
$j$ is sufficiently large, and so $F_j\in\rz$ and $f\sep F_j$. It is
well known that the zeros of a polynomial are continuous functions
of the coefficients of the polynomial and the limit of a sequence of
$\rz$ polynomials is still a $\rz$ polynomial (see \cite{Coo08} for
instance). Thus $F\in\rz$ and $f\sep F$ by continuity. Assume now
that $v(r)=0$ for some zero $r$ of $f$ with the multiplicity $m$.
Then $(x-r)^m|f$ implies $(x-r)^m|F$ from (\ref{Fff'}). Let
$f(x)=(x-r)^mh(x)$ and $F(x)=(x-r)^mH(x)$. Then $h(r)\neq 0$ and
$$H(x)=\left[u(x)+m\frac{v(x)}{x-r}\right]h(x)+v(x)h'(x)$$
by (\ref{Fff'}). So $H(r)=[u(r)+mv'(r)]h(r)$. If $u(r)+mv'(r)\neq
0$, then $H(r)\neq 0$, and so the multiplicity of $r$ as a zero of
$F$ is precisely $m$. This proves (b). If $u(r)+mv'(r)=0$, then
$H(r)=0$ and so the multiplicity of $r$ as a zero of $F$ is at least
$m+1$. However, $f\sep F$ and $r$ is a zero of $f$ with the
multiplicity $m$. Hence the multiplicity of $r$ as a zero of $F$ is
at most $m+1$. Thus the multiplicity of $r$ as a zero of $F$ is
precisely $m+1$. This proves (c).

(A)\quad Now let $r$ be the largest zero of $f$, with the
multiplicity $m$. Then $F$ has at most one zero larger than $r$
since $f\sep F$.

Assume that $v(r)\neq 0$. Then $r$ is a zero of $F$ with the
multiplicity $m-1$. Thus $F$ has one zero larger than $r$. Assume
that $v(r)=0$ and $u(r)+mv'(r)<0$. Then $h(r)>0$ since $h$ is
standard and has no zero larger than $r$. Hence
$H(r)=[u(r)+mv'(r)]h(r)<0$. Thus $H$ has one zero larger than $r$
since $H$ is standard, and so does $F$. Assume that $v(r)=0$ and
$u(r)+mv'(r)>0$. Then $H(r)>0$. Hence $H$ has an even number of
zeros larger than $r$. Thus $H$ has no zero larger than $r$, and so
does $F$. Assume that $v(r)=u(r)+mv'(r)=0$. Then $r$ is a zero of
$F$ with the multiplicity $m+1$. Thus $F$ has no zero larger than
$r$.

So we conclude that $F\in\rz(-\infty,r]$ if and only if $v(r)=0$ and
$u(r)+mv'(r)\ge 0$.

(B)\quad If $\deg F=\deg f$, then the result is clear since $f\sep
F$. If $\deg F=\deg f+1$, then let $g(x)=(-1)^nf(-x)$ and
$G(x)=(-1)^{n+1}F(-x)$ where $n=\deg f$. It follows that
$$G(x)=-u(-x)g(x)+v(-x)g'(x)$$
from (\ref{Fff'}). Thus the statement follows from (A).
\end{proof}

Combining (A) and (B) of Theorem~\ref{m-thm}, it is not difficult to
give a necessary and sufficient condition that guarantees zeros of
$f$ and $F$ are in the same closed interval. We omit the details for
the sake of brevity and only give the following result as a
demonstration.
\begin{corollary}\label{coro-x(x+1)}
Let $f$ and $F$ be two standard polynomials satisfying
\begin{equation*}
F(x)=(ax+b)f(x)+x(x+1)f'(x).
\end{equation*}
Suppose that $f(x)\in\rz[-1,0], x^{m_0}\|f$ and $(x+1)^{m_1}\|f$.
Then $b+m_0\ge 0$ and $a+m_1\ge b$ imply that $F\in\rz[-1,0]$ and
$f\sep F$. Furthermore, $x^{m_0}\|F$ if $b+m_0>0$ or $x^{m_0+1}\|F$
if $b+m_0=0$, and $(x+1)^{m_1}\|F$ if $a+m_1>b$ or
$(x+1)^{m_1+1}\|F$ if $a+m_1=b$.
\end{corollary}

Now we can apply Theorem~\ref{m-thm} to strengthen the results of
Brenti and Br\"and\'en about the $q$-Eulerian polynomials by the
recurrence (\ref{anxq-rr}) and by induction.
\begin{proposition}\label{e-prop}
Let $q\in\mathbb{R}$ and $n\in\mathbb{N}$.
\begin{enumerate}
\item [(a)]
If $q>0$, then $A_n(x;q)$ have nonpositive and simple zeros for
$n\ge 2$.
\item [(b)]
If $n+q\le 0$, then $A_{n+1}(x;q)\in\rz[1,+\infty)$.
\item [(c)]
If $q$ is a negative integer, then $A_{n}(x;q)\in\rz[1,+\infty)$ and
$(x-1)^{m}\|A_n(x;q)$ where $m=\max\{n+q,0\}$. In particular,
$A_n(x;-1)=-(x-1)^{n-1}$.
\end{enumerate}
\end{proposition}

We can also give an interpretation of the result when $q$ is a
negative integer. For this purpose, we give a $q$-analog of the
Frobenius formula of the classical Eulerian polynomials
$$A_n(x)=x\sum_{k=0}^{n}k!S(n,k)(x-1)^{n-k}.$$
\begin{proposition}
We have
\begin{equation}\label{q-Frobenius}
A_n(x;q)=\sum_{k=0}^{n}\binom{q+k-1}{k}k!S(n,k)(x-1)^{n-k}.
\end{equation}
\end{proposition}
\begin{proof}
We proceed  by induction on $n$. The equality is obvious for $n=0$
and $n=1$. Now assume that (\ref{q-Frobenius}) holds for $n\ge 1$.
Then by (\ref{anxq-rr}),
\begin{eqnarray*}
 A_{n+1}(x;q)
   &=& (nx+q)\sum_{k=0}^{n}\binom{q+k-1}{k}k!S(n,k)(x-1)^{n-k} \\
   & & +x(1-x)\sum_{k=0}^{n}\binom{q+k-1}{k}k!S(n,k)(n-k)(x-1)^{n-k-1}  \\
   &=& \sum_{k=0}^{n}\binom{q+k-1}{k}k!S(n,k)(kx+q)(x-1)^{n-k} \\
   &=& \sum_{k=0}^{n}\binom{q+k-1}{k}k!S(n,k)k(x-1)^{n-k+1}\\
   & & +\sum_{k=0}^{n}\binom{q+k-1}{k}k!S(n,k)(q+k)(x-1)^{n-k} \\
   &=&
   \sum_{k=0}^{n+1}\binom{q+k-1}{k}k!\left[kS(n,k)+S(n,k-1)\right](x-1)^{n-k+1}\\
   &=& \sum_{k=0}^{n+1}\binom{q+k-1}{k}k!S(n+1,k)(x-1)^{n-k+1}
\end{eqnarray*}
where we use the well-known recurrence $S(n+1,k)=kS(n,k)+S(n,k-1)$
for the Stirling numbers of the second kind in the last equality.
This completes the proof.
\end{proof}

When $q=-t$ is a negative integer, (\ref{q-Frobenius}) can be
written as
$$A_n(x;-t)=\sum_{k=0}^{n}(-1)^{k}\binom{t}{k}k!S(n,k)(x-1)^{n-k}.$$
It immediately follows that $(x-1)^{n-t}\|A_n(x;-t)$ for $n\ge t$,
as desired.
\section{Applications}
Theorem~\ref{m-thm} can provide a unified explanation of many known
results, including the fact that the classical Eulerian polynomials
and the Bell polynomials have only real simple zeros. In this
section we give more examples as applications.
\subsection{ Linear transformations preserving RZness}
Consider the invertible linear operator $\mathcal{T}:
\mathbb{R}[x]\rightarrow \mathbb{R}[x]$ defined by
$$\mathcal{T}{(x)_i}=x^i$$
for all $i\in\mathbb{N}$ and linear extension, where
$(x)_i=x(x-1)\cdots (x-i+1)$ and $(x)_0=1$. Wagner~\cite[Lemma
3.3]{Wag91} showed the following result.
\begin{proposition}\label{T-trans}
Let $\xi\in\mathbb{R}$ and $p$ be a real polynomial such that
$\mathcal{T}(p)\in\rz(-\infty,0]$. Then
\begin{enumerate}
  \item [\rm (a)] $F:=\mathcal{T}((x-\xi)p)\in\rz$.
  \item [\rm (b)] Let $m$ denote the multiplicity of $0$ as a
  zero of $\mathcal{T}(p)$. Then $F\in\rz(-\infty,0]$ if and only if $\xi\le m$.
  \item [\rm (c)] Furthermore, the multiplicity of $0$ as a zero
  of $F$ is $m$ if $\xi\neq m$, and is at least $m+1$ if $\xi=m$.
\end{enumerate}
\end{proposition}
Actually, let $f=\mathcal{T}(p)$. Then $F=(x-\xi)f+xf'$. Thus
Proposition~\ref{T-trans} is obvious from the viewpoint of
Theorem~\ref{m-thm}.

The {\it E-transformation} is the invertible linear operator
$\mathcal{E}: \mathbb{R}[x]\rightarrow \mathbb{R}[x]$ defined by
$$\mathcal{E}(\binom{x}{i})=x^i$$
for all $i\in\mathbb{N}$ and linear extension. This transformation
is important in the theory of $(P,\Omega)$-partitions (see
Brenti~\cite{Bre89} for details). Br\"and\'en~\cite[Lemma
4.4]{Bra06} showed the following.
\begin{proposition}\label{E-trans}
Let $\alpha\in[-1,0]$ and let $p$ be a polynomial such that
$\mathcal{E}(p)\in\rz[-1,0]$. Then
$\mathcal{E}((x-\alpha)p)\in\rz[-1,0]$ and $\mathcal{E}(p)\sep
\mathcal{E}((x-\alpha)p)$. If $\mathcal{E}(p)$ in addition only has
simple zeros, then so does $\mathcal{E}((x-\alpha)p)$.
\end{proposition}

Actually, let $f=\mathcal{E}(p)$ and $F=\mathcal{E}((x-\alpha)p)$.
Then $F=(x-\alpha)f+x(x+1)f'$. Thus Proposition~\ref{E-trans} is an
immediate consequence of Corollary~\ref{coro-x(x+1)}. Furthermore,
if the multiplicity of $0$ as a zero of $f$ is $m_0$, then the
multiplicity of $0$ as a zero of $F$ is also $m_0$ except
$m_0=\alpha=0$; if the multiplicity of $-1$ as a zero of $f$ is
$m_1$, then the multiplicity of $-1$ as a zero of $F$ is also $m_1$
except $m_1=0$ and $\alpha=-1$.
\subsection{Compositions of multisets}
Let $\bn=(n_1,n_2,\ldots)$ be the multiset consisting of $n_i$
copies of the $i$th type element. Denote by $\mo(\bn,k)$ the number
of compositions of $\bn$ into exactly $k$ parts. Then
\begin{equation}\label{onk-rr}
(n_j+1)\mo(\bn+e_j,k)=k\mo(\bn,k-1)+(n_j+k)\mo(\bn,k),
\end{equation}
where $\bn+e_j$ denotes the multiset obtained from $\bn$ by
adjoining one (additional) copy of the $j$th type element. Let
$f_{\bn}(x)=\sum_{k\ge 0}\mo(\bn,k)x^k$ be the associated generating
function. Then by (\ref{onk-rr}),
\begin{equation}\label{fnejx-rr}
(n_j+1)f_{\bn+e_j}(x)=(x+n_j)f_{\bn}(x)+x(x+1)f'_{\bn}(x).
\end{equation}
Simion showed that the multiplicity of $-1$ as a zero of
$f_{\bf}(x)$ is $\max_i\{n_i-1\}$ by means of the theory of posets
(\cite[Lemma 1.1]{Sim84}). Based on this result and appropriate
transformation to the recurrence (\ref{fnejx-rr}), she further
showed that $f_{\bn}(x)\in\rz[-1,0]$ and $f_{\bn}(x)\sep
f_{\bn+e_j}(x)$ (\cite[Theorem 1]{Sim84}). These results are now
clear from the viewpoint of Corollary~\ref{coro-x(x+1)}.

In particular, if $\bn=(1,1,\ldots,1)$, then $\mo(\bn,k)=k!S(n,k)$,
where $S(n,k)$ is the Stirling number of the second kind. Thus the
polynomial $F_n(x)=\sum_{k=1}^{n}k!S(n,k)x^k$ has only real simple
zeros in the interval $[-1,0]$. It is interesting that
$F_n(x)=\frac{x^{n+1}}{x+1}A_n(\frac{x+1}{x})$ by the Frobenius
formula, where $A_n(x)$ is the classical Eulerian polynomial.
\subsection{Alternating runs}
Let $\pi=a_1a_2\cdots a_n\in\sn$. We say that $\pi$ changes
direction at position $i$ if either $a_{i-1}<a_i>a_{i+1}$, or
$a_{i-1}>a_i<a_{i+1}$. We say that $\pi$ has $k$ {\it alternating
runs} if there are $k-1$ indices $i$ such that $\pi$ changes
direction at these positions. Let $R(n,k)$ denote the number of
permutations in $\sn$ having $k$ alternating runs. Then
\begin{equation}\label{rnk-rr}
    R(n,k)=kR(n-1,k)+2R(n-1,k-1)+(n-k)R(n-1,k-2)
\end{equation}
for $n,k\ge 1$, where $R(1,0)=1$ and $R(1,k)=0$ for $k\ge 1$ (see
B\'ona~\cite[Lemma 1.37]{Bona04} for a combinatorial proof). Let
$R_n(x)=\sum_{k=1}^{n-1}R(n,k)x^k$.
Then the recurrence (\ref{rnk-rr}) induces
\begin{equation}\label{rnx-rr}
R_{n+2}(x)=x(nx+2)R_{n+1}(x)+x\left(1-x^2\right)R_{n+1}'(x),
\end{equation}
with $R_1(x)=1$ and $R_2(x)=2x$. B\'ona and Ehrenborg~\cite[Lemma
2.3]{BE00} showed that $R_n(x)$ has the zero $x=-1$ with
multiplicity $\lrf{\frac{n}{2}}-1$ and suspected that the other half
zeros of $R_n(x)$ are all real, negative and distinct. The
polynomial $R_n(x)$ is closely related to the classical Eulerian
polynomial $A_n(x)$:
\begin{equation}\label{rnx-anx}
R_n(x)=\left(\dfrac{1+x}{2}\right)^{n-1}(1+w)^{n+1}A_n\left(\dfrac{1-w}{1+w}\right),\quad
w=\sqrt{\frac{1-x}{1+x}}
\end{equation}
(Knuth~\cite[p. 605]{Knuth73}). From the relation (\ref{rnx-anx})
and the fact that $A_n(x)$ have only real zeros, Wilf can show that
$R_n(x)$ have only real zeros for $n\ge 2$ (see B\'ona~\cite[Theorem
1.41]{Bona04} and Stanley~\cite{Sta05} for details). Very recently,
Canfield and Wilf~\cite{CW0609704} pointed out (without proof) that
this result can also be obtained based on the recurrence
(\ref{rnx-rr}). Indeed, we can give the following more precise
result by Theorem~\ref{m-thm}.
\begin{proposition}\label{Bona-Wilf}
Let $R_n(x)$ be the generating function of alternating runs. Then
$R_n(x)\in\rz[-1,0]$ and $R_n(x)\sep R_{n+1}(x)$ for $n\ge 1$. More
precisely, $R_n(x)$ has $\lrc{\frac{n}{2}}$ simple zeros including
$x=0$, and the zero $x=-1$ with the multiplicity
$\lrf{\frac{n}{2}}-1$.
\end{proposition}

\end{document}